\documentclass[12pt]{amsart}
\usepackage{amsmath}
\usepackage{amssymb}
\usepackage{amsmath,amssymb,amsthm,amscd}
\usepackage{relsize}

\usepackage{txfonts}
\usepackage{amsbsy}
\usepackage{graphicx,color,epstopdf}

\numberwithin{equation}{section}

\newtheorem{prop}{Proposition}[section]
\newtheorem{theo}{Theorem}[section]
\newtheorem{lemm}{Lemma}[section]

\def\begeq{\begin{equation}}
\def\endeq{\end{equation}}

\begin{document}

\title{Necessary and Sufficient Conditions to Bernstein Theorem of a Hessian Equation}
\author{Shi-Zhong Du}
\thanks{The author is partially supported by NSFC (12171299), and GDNSF (2019A1515010605)}
  \address{The Department of Mathematics,
            Shantou University, Shantou, 515063, P. R. China.} \email{szdu@stu.edu.cn}

\renewcommand{\subjclassname}{%
  \textup{2010} Mathematics Subject Classification}
\date{Feb. 2021}

\begin{abstract}
   The Hessian quotient equations
      \begin{equation}\label{e0.1}
        S_{k,l}(D^2u)\equiv\frac{S_k(D^2u)}{S_l(D^2u)}=1, \ \ \forall x\in{\mathbb{R}}^n
      \end{equation}
   were studied for $k-$th symmetric elementary function $S_k(D^2u)$ of eigenvalues $\lambda(D^2u)$ of the Hessian matrix $D^2u$, where $0\leq l<k\leq n$. For $l=0$, \eqref{e0.1} is reduced to a $k-$Hessian equation
      \begin{equation}\label{e0.2}
         S_k(D^2u)=1, \ \ \forall x\in{\mathbb{R}}^n.
      \end{equation}
   Two quadratic growth conditions were found by Bao-Cheng-Guan-Ji (\cite{BCGJ}, American J. Math., 2003, {\bf125}, 301-316) ensuring the Bernstein properties of \eqref{e0.1} and \eqref{e0.2} respectively. In this paper, we will drop the point wise quadratic growth condition of \cite{BCGJ} and prove three necessary and sufficient conditions to Bernstein property of \eqref{e0.1} and \eqref{e0.2}, using a reverse isoperimetric type inequality, volume growth or $L^p$-integrability respectively. Our new volume growth or $L^p-$integrable conditions improve largely various previously known point wise conditions in \cite{BCGJ,CX,CY,LRW,Y} etc.
\end{abstract}

\maketitle\markboth{Bernstein problem}{$\gamma-$ball condition}

\tableofcontents

\section{Introduction}

Fully nonlinear equation in form of
  \begin{equation}\label{e1.1}
      F(D^2u)=\psi, \ \ \forall x\in\Omega
  \end{equation}
has been studied extensively in the past years, where
    $$
     F(D^2u)=f(\lambda_1,\cdots,\lambda_n)
    $$
is assumed to be a symmetric function on eigenvalues
   $$
    \lambda(D^2u)\equiv(\lambda_1,\cdots,\lambda_n)\in{\mathbb{R}}^n.
   $$
It is well known that \eqref{e1.1} is elliptic if
   $$
     \frac{\partial f}{\partial\lambda_i}>0, \ \ \forall i=1,2,\cdots,n.
   $$
Moreover, $F(D^2u)$ is a concave function on $D^2u\in{\mathbb{R}}^{n^2}$ if $f(\lambda)$ is a concave function on $\lambda$ \cite{CNS}. When studying \eqref{e1.1}, a natural class of solutions ensuring the ellipticity of the equation is the so called $k-$admissible function defined as following. We say a function $u\in C^2(\Omega)$ is $k$-admissible if
   \begin{equation}\label{e1.2}
     \lambda(D^2u)\in\Gamma_k,
   \end{equation}
where $\Gamma_k$ is the G{\aa}rding convex cone on ${\mathbb{R}}^n$ defined by
  \begin{equation}\label{e1.3}
   \Gamma_k\equiv\Big\{\lambda=(\lambda_1,\cdots,\lambda_n)\in{\mathbb{R}}^n\Big|\ S_j(\lambda)>0, \ \ \forall j=1,\cdots, k\Big\}
  \end{equation}
and
  $$
   S_j(\lambda)\equiv\Sigma_{1\leq i_1<i_2\cdots<i_j\leq n}\lambda_{i_1}\lambda_{i_2}\cdots\lambda_{i_j}
  $$
is the $j-$th elementary symmetric polynomial of $\lambda$ for each $j=1,2,\cdots,n$. As usually, we denote $S_0(D^2u)\equiv1$ for simplicity. It is clear that
  $$
   \Gamma_1=\Big\{\lambda\in{\mathbb{R}}^n\Big|\ \Sigma_{j=1}^n\lambda_j>0\Big\}
  $$
is the half space and
  $$
   \Gamma_n=\Big\{\lambda\in{\mathbb{R}}^n\Big|\ \lambda_j>0, \ \ \forall j=1,\cdots,n\Big\}
  $$
corresponds exactly to the family of convex functions. In generally, the cone $\Gamma_k$ can also be equivalently defined as the component of
  $$
   \Big\{\lambda\in{\mathbb{R}}^n\Big|\ S_k(\lambda)>0\Big\}
  $$
containing the vector $(1,\cdots,1)$, or be characterized as
  \begin{equation}\label{e1.4}
   \Gamma_k=\Big\{\lambda\in{\mathbb{R}}^n\Big| 0<S_k(\lambda)\leq S_k(\lambda+\eta), \ \ \forall \eta_j\geq0, j=1,\cdots,n\Big\}.
  \end{equation}
Supposing that $D^2u$ is diagonal at some point and its eigenvalues are different from each other, there holds
  \begin{equation}\label{e1.5}
    \frac{\partial\lambda_i}{\partial u_{pq}}=\delta_{pi}\delta_{qi},\ \ \frac{\partial^2\lambda_i}{\partial u_{pq}\partial u_{rs}}=\frac{1}{\lambda_i-\lambda_q}\delta_{pi}\delta_{ri}\delta_{sq}(1-\delta_{pq}).
  \end{equation}
As a result, one gets that
  \begin{eqnarray}\label{e1.6}
    &\displaystyle\frac{\partial F}{\partial u_{pq}}=\frac{\partial f}{\partial\lambda_p}\delta_{pq},& \\ \nonumber
     &\displaystyle\frac{\partial^2F}{\partial u_{pq}\partial u_{rs}}=\frac{f_p-f_q}{\lambda_p-\lambda_q}\delta_{pr}\delta_{qs}(1-\delta_{pq})+\frac{\partial^2f}{\partial\lambda_p\partial\lambda_r}\delta_{pq}\delta_{rs}&
  \end{eqnarray}
for diagonal matrix $D^2u$.

In this paper, we study the Hessian quotient equation of the form
     \begin{equation}\label{e1.7}
       S_{k,l}(D^2u)\equiv\frac{S_k(D^2u)}{S_l(D^2u)}=1, \ \ \forall x\in{\mathbb{R}}^n
     \end{equation}
for $0\leq l<k\leq n$. If $l=0$ and $k=n$, \eqref{e1.7} becomes the classical Monge-Amp\`{e}re equation
   \begin{equation}\label{e1.8}
     \det(D^2u)=1, \ \ \forall x\in{\mathbb{R}}^n.
   \end{equation}
A well known Bernstein property asserts that any locally strict convex smooth solution of \eqref{e1.8} must be a quadratic function. It was shown by J\"{o}rgens \cite{J} for $n=2$, Calabi \cite{C} for $n=3,4,5$ and Pogorelov \cite{P} for all dimensions $n\geq2$. The result was later extended to viscosity solutions by Caffarelli \cite{Ca}. The readers may also refer to Cheng-Yau \cite{CY} for another geometry proof. If $l=0, k=1$, \eqref{e1.7} is exactly the linear Poisson's equation
   \begin{equation}\label{e1.9}
     \triangle u=1, \ \ \forall x\in{\mathbb{R}}^n.
   \end{equation}
When considering convex solutions of \eqref{e1.9}, one can still verify the validity of Bernstein property without difficulty. In fact, given any unit vector $\xi\in{\mathbb{S}}^{n-1}$, by convexity of $u$, the harmonic function $u_{\xi\xi}$ is bounded from below by constant zero. After applying the Liouville theorem to nonnegative harmonic function, one concludes that $u$ must be a quadratic function.

When $l=0, 2\leq k\leq n-1$, \eqref{e1.7} changes to the $k$-Hessian equation
     \begin{equation}\label{e1.10}
       S_k(D^2u)=1, \ \ \forall x\in{\mathbb{R}}^n.
     \end{equation}
Unlike the Monge-Amp\`{e}re equation \eqref{e1.8} or Poisson equation \eqref{e1.9}, the Bernstein problem of \eqref{e1.10} for $2\leq k\leq n-1$ is much more complicate to be explored. A first sufficient condition was found by Bao-Chen-Guan-Ji \cite{BCGJ} using point wise quadratic growth form.

\begin{theo}\label{t1.1}
  Let $u\in C^\infty({\mathbb{R}}^n)$ be a locally strict convex solution of the $k-$Hessian equation \eqref{e1.10} for $2\leq k\leq n-1$. Supposing that
    \begin{equation}\label{e1.11}
      u(x)\geq a_1|x|^2-a_2, \ \ \forall x\in{\mathbb{R}}^n
    \end{equation}
  holds for some positive constants $a_1, a_2$, then $u$ must be a quadratic polynomial.
\end{theo}

The result was later extended by Chang-Yuan \cite{CYu} to uniform convexity assumption
   $$
    D^2u\geq\Bigg(\delta-\sqrt{\frac{2}{n(n-1)}}\Bigg)I, \ \ \forall x\in{\mathbb{R}}^n
   $$
for $k=2$ and some positive constant $\delta$. At the same time, they guess their result should still hold true under semiconvexity assumption
   $$
    D^2u(x)\geq -KI, \ \ \forall x\in{\mathbb{R}}^n
   $$
for arbitrarily large $K$, even for $2\leq k\leq n-1$. The validity of this conjecture was verified later by Yuan \cite{Y} for $n=3$ and $k=2$. More recently, Li-Ren-Wang \cite{LRW} obtained a Bernstein result of \eqref{e1.10} under the assumptions of $(k+1)-$admissibility and quadratic growth. If restricting to the case $k=2$, Chen-Xiang \cite{CX} have shown a Bernstein property using the conditions of $2-$admissible, quadratic growth and
   $$
    S_3(D^2u(x))\geq -A, \ \ \forall x\in{\mathbb{R}}^n.
   $$
The last condition of lower bound of sigma-3 can be removed when $n=3$.

 Turning to the Hessian quotient equation \eqref{e1.7}, Bao-Chen-Guan-Ji have also obtained the following Bernstein property in \cite{BCGJ}.

\begin{theo}\label{t1.2}
  Let $u\in C^\infty({\mathbb{R}}^n)$ be a locally strict convex solution of \eqref{e1.7} for $1\leq l<k=n$. If
    \begin{equation}\label{e1.12}
      u(x)\leq a_3(|x|^2+1), \ \ \forall x\in{\mathbb{R}}^n
    \end{equation}
  holds for some positive constant $a_3$, then $u$ must be a quadratic function. Moreover, if $l=n-1$, Theorem \ref{t1.2} holds without assumption \eqref{e1.12}.
\end{theo}

At the same paper, the authors laid down two open problems about the removability of growth assumptions in Theorem \ref{t1.1} and \ref{t1.2}. (see \cite{BCGJ}, Page 302, Remark 1.3 for \eqref{e1.7}; and Page 314 for \eqref{e1.10})\  The first main purpose of this paper is to remove the condition of point wise quadratic growth and prove the following necessary and sufficient condition to Bernstein problem of \eqref{e1.10}.

\begin{theo}\label{t1.3}
  Considering the $k-$Hessian equation \eqref{e1.10} for $2\leq k\leq n-1$, any locally strict convex solution $u\in C^\infty({\mathbb{R}}^n)$ satisfying the normalized condition
    \begin{equation}\label{e1.13}
      u(0)=0, \ \ Du(0)=0
    \end{equation}
  is a quadratic function if and only if a reverse isoperimetric inequality
    \begin{equation}\label{e1.14}
      {\mathcal{I}}(u)\equiv\liminf_{t\to\infty}\frac{\int_{\Omega_t}|D u|}{\Big(\int_{\Omega_t}|t-u|^{\frac{n}{n-1}}\Big)^{\frac{n-1}{n}}}<\infty
    \end{equation}
  is satisfied, where
     $$
      \Omega_t\equiv\Big\{x\in{\mathbb{R}}^n\Big|\ u(x)<t\Big\}
     $$
  is the sub-level set of the function $u$.
\end{theo}

It is worthy to mention that the isoperimetric inequality
   $$
    \Bigg(\int_{\Omega_t}|t-u|^{\frac{n}{n-1}}\Bigg)^{\frac{n-1}{n}}\leq C_n\int_{\Omega_t}|D u|
   $$
was well known for bounded variation functions. Moreover, our necessary and sufficient condition of reverse isoperimetric inequality \eqref{e1.14} is not a restrictive one due to the following gradient inequality by Trudinger-Wang \cite{TW}.

\begin{theo}\label{t1.4}
  Let $u\in C^2(\Omega)\cap\Gamma_k, k=1,\cdots,n$ satisfying $u\leq 0$ in $\Omega$. Then for any subdomain $\Omega'\subset\subset\Omega$, there hold the estimates
     \begin{equation}\label{e1.15}
       \int_{\Omega'}|Du|^qS_l(D^2u)\leq Cd^{-2l-q}_{\Omega'}\Bigg(\int_\Omega|u|\Bigg)^{q+l}, \ \ d_{\Omega'}\equiv dist(\Omega',\partial\Omega)
     \end{equation}
  for all $l=0,\cdots,k-1, 0\leq q<\frac{n(k-l)}{n-k}$, where $C$ is a constant depending only on $n, k, l$ and $q$.
\end{theo}

Our second purpose is to prove a necessary and sufficient condition of Bernstein problem of Hessian quotient equation \eqref{e1.7} in case of $1\leq l\leq n-2$ and $k=n$. In fact, we first reduce the Bernstein problem to the following interior $C^2-$estimation of $(n-l)-$Hessian equation without strict convexity, which is comparable to the version of $C^2-$theorem of Chou-Wang in \cite{CW}.

\begin{theo}\label{t1.5}
     Supposing that there exists some positive function $\psi:\ {\mathbb{R}}^3\to{\mathbb{R}}^+$, such that for any convex solution $v$ of
     \begin{equation}\label{e1.16}
      \begin{cases}
         S_{n-l}(D^2v)=1, & \forall y\in B_1(0),\\
         v(0)=0, Dv(0)=0,
      \end{cases}
     \end{equation}
there holds
    \begin{equation}\label{e1.17}
       |D^2v|(0)\leq \psi(n,l,\sup_{B_1}|Dv|),
    \end{equation}
then any locally strict convex solution $u$ of \eqref{e1.7} for $1\leq l<k=n$ must be a quadratic function if and only if reverse isoperimetric inequality \eqref{e1.14} holds upon normalized condition \eqref{e1.13}.
\end{theo}

The assumption of Theorem \ref{t1.5} is clearly true for $l=n-1$. To handle the case $l=n-2$, we need the following striking $C^2$-estimation, which was due to Warren-Yuan \cite{WY} for $n=2$, Mc-Song-Yuan \cite{MSY} for convex solutions and all dimensions $n\geq2$.

\begin{theo}\label{t1.6}
  Let $u$ be a smooth convex solution to \eqref{e1.16} on unit ball $B_1\subset{\mathbb{R}}^n$ for $l=n-2$. Then we have
    \begin{eqnarray}\nonumber\label{e1.18}
      \sup_{B_{1/4}}|D^2u|&\leq& C_n\exp\Big(C_n\max_{B_{1/2}}|Du|^3\Big)\\
        &\leq&C_n\exp\Big(C_n\max_{B_{1}}|u|^3\Big)
    \end{eqnarray}
  for some universal constant $C_n$.
\end{theo}

The result was later generalized to semi-convex solutions by Shankar-Yuan \cite{SY}, and generalized to equations
   \begin{equation}\label{e1.19}
     S_2(D^2u)=f(x,u,Du), \ \ \forall x\in\Omega
   \end{equation}
by Guan-Qiu \cite{GQ}. As a corollary of Theorems \ref{t1.5} and \ref{t1.6},  we derive the following result.

\begin{theo}\label{t1.7}
  Considering the Hessian quotient equation \eqref{e1.7} for $k=n$, $l=n-1$ or $l=n-2$, any smooth locally strict convex solution $u$ must be a quadratic polynomial if and only if reverse isoperimetric inequality \eqref{e1.14} holds upon normalized condition \eqref{e1.13}.
\end{theo}

It's remarkable that for $k=n=3$ and $l=1$, the Hessian quotient equation \eqref{e1.7} is equivalent to the special Lagrangian equation
   $$
    \arctan\lambda_1+\arctan\lambda_2+\arctan\lambda_3=\pi.
   $$
So, after applying an interesting result of Yuan in \cite{Y2}, Bernstein property in Theorem \ref{t1.7} can be shown without the help of reverse isoperimetric inequality.

 Since the reverse isoperimetric inequality is not easy to be verified, in the final part of this paper, we shall also reduce \eqref{e1.14} to equivalent conditions of volume growth or $L^p$-integrability in the following sense.

\begin{theo}\label{t1.8}
  Under the assumptions of Theorems \ref{t1.3}, \ref{t1.5} or \ref{t1.7} and normalized condition \eqref{e1.13}, the reverse isoperimetric inequality \eqref{e1.14} is equivalent to the following two alternative conditions
    \begin{equation}\label{e1.20}
      \liminf_{t\to\infty}t^{-n/2}|\Omega_t|<\infty
    \end{equation}
  or
    \begin{equation}\label{e1.21}
      \liminf_{t\to\infty}t^{-p-n/2}\int_{\Omega_t}(u+1)^p<\infty
    \end{equation}
  for some $p\in{\mathbb{R}}$, where $|\Omega_t|$ is the volume of sub-level set of $u$.
\end{theo}

Our necessary and sufficient condition \eqref{e1.20} or condition \eqref{e1.21} improve largely various known point wise conditions in \cite{BCGJ,CX,CY,LRW,Y} etc. As the $k-$admissible structure of Hessian quotient equation \eqref{e1.7}, it is natural to ask whether any $k-$admissible solution of \eqref{e1.7} for $0\leq l<k\leq n-1$ must be a quadratic polynomial. A surprising example shown by Warren \cite{W} gives non-polynomial entire $k$-admissible solutions in case $l=0$ and $1\leq k\leq\frac{n+1}{2}$, for which reverse isoperimetric inequality \eqref{e1.14}, volume growth condition \eqref{e1.20} or $L^p-$integrable condition \eqref{e1.21} all failed to hold. Therefore, one may not expect a general confirmed answer to this question. It would be interesting to ask whether there exist some non-quadratic convex solutions of \eqref{e1.7} such that reverse isoperimetric inequality, volume growth or $L^p-$integrability fail to hold.

Our contents are organized as follows. We will first introduce a $\gamma-$ball condition and then prove some useful geometric lemmas in Section 2. With the help of a radius estimation, we complete the proof of Theorem \ref{t1.3} in Section 3. Next, we will give the proof to Bernstein Theorem \ref{t1.5} and its corollary Theorem \ref{t1.7} in Section 4. In the Section 5, we prove a key relation between reverse isoperimetric inequality and $\gamma-$ball condition in Theorem \ref{t5.1}. Finally, the equivalence of reverse isoperimetric inequality with volume growth or $L^p-$integrable conditions will be presented in the last Section 6.

\vspace{40pt}

\section{$\gamma-$ball condition and several crucial geometric lemmas}

Letting $u$ be a convex solution of \eqref{e1.7} or \eqref{e1.10}, after subtracting a linear function, one may always assume that the normalized condition
   $$
    u(0)=0, \ \ Du(0)=0
   $$
holds throughout this paper. To proceed further, for each given $t_0>0$, let us normalize $u$ to be
   $$
     u^a(x)\equiv t_0^{-1}u(\sqrt{t_0}x)
   $$
in spirit of Pogorelov \cite{P} and set
   $$
    \Omega^a_t\equiv\Big\{x\in{\mathbb{R}}^n\Big|\ u^a(x)<t\Big\}, \ \ t>0
   $$
to be its sub-level sets. Original normalizations of Pogorelov contain a second normalization
    $$
     u^b(x)\equiv u^a(A^{-1}x)
    $$
by affine transformation $y=Ax, \det A=1$. However, $k-$Hessian polynomial is not invariant under the affine transformations in general. So, the John's lemma can not be applied directly. Fortunately, for any $\gamma>1$, one can still define a $\gamma-$ball condition to $\Omega^a\equiv\Omega^a_1$ by
  \begin{equation}\label{e2.1}
   \gamma^{-1}B_R\subset\Omega^a\subset \gamma B_R, \ \ B_R\equiv B_R(x_0).
  \end{equation}
A bounded convex domain satisfying the $\gamma-$ball condition \eqref{e2.1} is called to be a $\gamma-$ball domain for short. By John's lemma \cite{Jo}, there exist an universal constant $C_n$ and an affine transformation $A, \det(A)=1$ with eigenvalues
   $$
    \mu(A)=(\mu_1,\mu_2,\cdots,\mu_n), \ \ 0<\mu_1\leq\mu_2\leq\cdots\leq\mu_n,
   $$
such that
   \begin{equation}\label{e2.2}
     C_n^{-1}B_R\subset A(\Omega^a)\subset C_n B_R.
   \end{equation}
It is clear that
   $$
    C_n^{-1}\gamma^2\leq \mu_n/\mu_1\leq C_n\gamma^2
   $$
holds for another universal constant $C_n$. Conversely, if $\Omega^a$ satisfies the John's ball condition for some affine transformation $A$ with eigenvalues $\mu(A)$, then $\Omega^a$ satisfies $\gamma-$ball condition for
   $$
    C_n^{-1}\sqrt{\frac{\mu_n}{\mu_1}}\leq\gamma\leq C_n\sqrt{\frac{\mu_n}{\mu_1}}.
   $$
In order showing our geometric lemmas, one needs to define the following normal mapping
   $$
    \nu(x)\equiv(D_1u,D_2u,\cdots, D_nu,-1), \ \ \forall x\in\Omega
   $$
from a given convex function $u$ as usually, and use the fact that the area of the image of normal mapping of $u$ equals to
   \begin{equation}\label{e2.3}
     \nu(u,\Omega)\equiv\int_{\nu(\Omega)}d\nu=\int_\Omega\det D^2u dx.
   \end{equation}
Now, let us prove several crucial geometric lemmas which would be useful later.

 \begin{lemm}\label{l2.1}
   Letting $u^a\in C^2(\Omega^a)$ be the convex function after Pogorelov's first normalization,
     \begin{equation}\label{e2.4}
       R^{-n}\leq C_{n,1}\gamma^{n}\int_{\Omega^a}\det D^2u^adx
     \end{equation}
   holds for positive constant $C_{n,1}\equiv\omega_n^{-1}$, where $\omega_n\equiv|B_1|$ is the volume of $n$-dimensional unit ball.
 \end{lemm}

 \noindent\textbf{Proof.} At first, we construct a cone $V_1$ carrying the origin as its vertex and carrying
   $$
     \widetilde{\Omega}^a\equiv\Big\{(x,u)\in{\mathbb{R}}^{n+1}\big|\ x\in\Omega^a, u=1\Big\}
   $$
 as its base. By convexity of $u^a$, the normal image of $V_1$ must be contained inside of the normal image of graph
    $$
     G\equiv\Big\{(x,u)\in{\mathbb{R}}^{n+1}\big|\ x\in\Omega^a, u=u^a(x)\Big\}.
    $$
 Another hand, since $\Omega^a$ is contained inside a large ball $B_{\gamma R}$, the area of the normal image of $V_1$ is no less than $\omega_n(\gamma R)^{-n}$. So, we achieve \eqref{e2.4} by \eqref{e2.3} and comparison. $\Box$\\

 The second lemma estimates $R$ from above.

 \begin{lemm}\label{l2.2}
    Letting $u^a\in C^2(\Omega^a)$ be the convex function after Pogorelov's first normalization,
     \begin{equation}\label{e2.5}
       R^{-n}\geq C_{n,2}^{-1}\gamma^{-n}\int_{B_{R/(2\gamma)}(x_0)}\det D^2u^adx
     \end{equation}
    holds for positive constant $C_{n,2}\equiv 2^n\omega_n$.
 \end{lemm}

 \noindent\textbf{Proof.} Let us take $(x_0,-1)$ to be the vertex and take $\widetilde{\Omega}^a$ defined in proof of Lemma \ref{l2.1} to be the base of a cone $V_2$. Noting that the interior of $V_2$ contains a portion $G_0$ of the graph $G$ of $u^a$, whose projection
    $$
     \Omega^a_0\equiv\Big\{x\in\Omega\big|\ (x,u^a(x))\in G_0\Big\}
    $$
 contains $\Omega^a/2\supset B_{R/(2\gamma)}(x_0)$, one concludes that \eqref{e2.5} holds by comparison of normal image. $\Box$\\

 A third lemma gives $C^1$-estimation of $u^a$.

 \begin{lemm}\label{l2.3}
   Letting $u^a\in C^2(\Omega^a)$ be the convex function after Pogorelov's first normalization, there exists an universal constant $C_{n,3}>0$ such that for any $0<s<t\leq 1$, there holds
     \begin{equation}\label{e2.6}
       \sup_{\Omega^a_s}|Du^a|\leq C_{n,3}\Bigg(\frac{\gamma R}{t-s}\Bigg)^{n-1}\int_{\Omega^a_t}\det D^2u^adx,
     \end{equation}
   where $\Omega^a_t\equiv\Big\{x\in{\mathbb{R}}^n\big|\ u^a(x)<t\Big\}$.
 \end{lemm}

 \noindent\textbf{Proof.} For any $x\in\Omega^a_s$, let
   $$
     h\equiv t-u^a(x)\geq t-s, \ \ d\equiv dist(x,\partial\Omega^a_t)
   $$
 and $D\leq2\gamma R$ be the diameter of $\Omega^a_t$. It is clear that $|Du^a|(x)\leq\frac{h}{d}$. Constructing a cone $V_3$ with vertex $(x,u^a(x))$ and base $\Omega^a_t$, the normal image of $V_3$ is contained inside of the normal image of
    $$
     G^a_t\equiv\Big\{(x,u^a(x))\big|\ x\in\Omega^a_t\Big\}.
    $$
 Because the normal image of a convex function is a convex set, and the normal image of $V_3$ contains a point $P$ of distant $h/d$ from origin and a ball $B_{h/D}(0)$ of radius $h/D$, we conclude that the cone with vertex $P$ and bottom $B_{h/D}(0)$ is contained inside of the normal image of $V_3$. So, one derives that
   $$
    C_{n,3}^{-1}\Bigg(\frac{h}{d}\Bigg)\Bigg(\frac{h}{D}\Bigg)^{n-1}\leq\int_{\Omega^a_t}\det D^2u^adx
   $$
 and hence \eqref{e2.6}, for $h\geq t-s$ and $D\leq 2\gamma R$. $\Box$\\

\vspace{40pt}

\section{Radius estimation and the proof of Theorem \ref{t1.3}}

When considering $k-$Hessian equation \eqref{e1.10}, Lemmas \ref{l2.1}-\ref{l2.3} are not sufficient to produce the desired $C^1$-bound of the solution since $R$ may be large in Lemma \ref{l2.3}. Fortunately, with the help of convexity and the equation \eqref{e1.10}, one can still estimate the radius $R$ from above.

\begin{prop}\label{p3.1}
  Suppose that the Dirichlet problem
    \begin{equation}\label{e3.1}
      \begin{cases}
        S_k(D^2u)=1, & \forall x\in\Omega\ni\{0\}\\
        u(0)=0, Du(0)=0, u(x)=1, & \forall x\in\partial\Omega
      \end{cases}
    \end{equation}
  admits a smooth convex solution $u$ for some convex domain $\Omega$ satisfying $\gamma-$ball condition
     \begin{equation}\label{e3.2}
      \gamma^{-1}B_R\subset\Omega\subset \gamma B_R, \ \ B_R\equiv B_R(x_0).
     \end{equation}
  Then there exists an universal constant $C_{n,4}>0$, such that
    \begin{equation}\label{e3.3}
       R\leq C_{n,4}\gamma.
    \end{equation}
\end{prop}

\noindent\textbf{Proof.} By Newton-Maclaurin's inequality, we have
   $$
    \frac{\triangle u}{n}\geq\Bigg(\frac{S_k(D^2u)}{C_n^k}\Bigg)^{\frac{1}{k}}=(C_n^k)^{-\frac{1}{k}}.
   $$
Therefore, $u$ must be a sub-solution of
   \begin{equation}\label{e3.4}
     \begin{cases}
       \triangle u\geq\alpha_{n,k}\equiv n(C_n^k)^{-\frac{1}{k}}, & \forall x\in\Omega\\
       u(x)=1, & \forall x\in\partial\Omega.
     \end{cases}
   \end{equation}
Noting that
   $$
    v(x)\equiv\frac{\alpha_{n,k}}{2n}\Big(|x-x_0|^2-\gamma^{-2}R^2\Big)+1
   $$
is a solution to
   $$
     \triangle v=\alpha_{n,k}, \ \ \forall x\in \gamma^{-1}B_R
   $$
satisfying
    $$
     v(x)=1, \ \ \forall x\in\partial\Big(\gamma^{-1}B_R\Big).
    $$
Comparing $u$ with $v$ by the maximum principle, one obtains that
   \begin{equation}\label{e3.5}
     u(x)\leq v(x)=\frac{\alpha_{n,k}}{2n}\Big(|x-x_0|^2-\gamma^{-2}R^2\Big)+1, \ \ \forall x\in \gamma^{-1}B_R.
   \end{equation}
Taking $x=x_0$ in \eqref{e3.5}, it yields that
   $$
    0\leq u(x_0)\leq1-\frac{\alpha_{n,k}}{2n \gamma^2}R^2.
   $$
Setting
   $$
    C_{n,4}\equiv\sqrt{\frac{2n}{\alpha_{n,k}}},
   $$
the proof of \eqref{e3.3} was done. $\Box$\\

 Now, let us turn back to prove Theorem \ref{t1.3}. We need first to quote a theorem of $C^2-$estimation by Chou-Wang \cite{CW}.

\begin{theo}\label{t3.1}
  Letting $u$ be a smooth convex solution of \eqref{e3.1}-\eqref{e3.2}, there exists a positive constant $C$ depending only on $n, k$ and $||u||_{C^1(\Omega)}$ but not on $\Omega$, such that
    \begin{equation}\label{e3.6}
      (u-1)^4|D^2u(x)|\leq C, \ \ \forall x\in\Omega.
    \end{equation}
\end{theo}

\noindent\textbf{Continue to prove Theorem \ref{t1.3}.} Utilizing the Newton-Maclaurin's inequality again, one gets that
   \begin{equation}\label{e3.7}
       \det D^2u^a\leq\Bigg(\frac{S_k(D^2u^a)}{C_n^k}\Bigg)^{\frac{n}{k}}=(C_n^k)^{-\frac{n}{k}}\equiv\beta_{n,k}, \ \ \forall x\in\Omega^a.
   \end{equation}
Combining Proposition \ref{p3.1} with Lemma \ref{l2.3}, it is inferred from \eqref{e3.7} that
   \begin{equation}\label{e3.8}
     |Du^a|\leq C_{1/2}, \ \ \forall x\in\Omega^a_{1/2},
   \end{equation}
where $C_{1/2}$ is a positive constant depending on $n,k,\gamma$ but not on $u^a$ and $\Omega^a$. Using Theorem \ref{t3.1} for $k-$Hessian equation \eqref{e1.10}, we conclude that
   \begin{equation}\label{e3.9}
     ||u^a||_{C^2(\Omega^a_{1/4})}\leq C_{1/4}, \ \ \forall x\in\Omega^a_{1/4}.
   \end{equation}
Once the second order derivative estimation \eqref{e3.9} has been obtained, \eqref{e1.10} becomes a uniformly elliptic equation. By imposing the Krylov's regularity theory \cite{K}, we reach the following higher regularity of the normalized solution.

\begin{prop}\label{p3.2}
  Let $u=u^a\in C^\infty(\Omega^a)$ be a strict convex solution of \eqref{e1.10} after Pogorelov's first normalization. Supposing that $\Omega^a$ satisfies the $\gamma-$ball condition, then
    \begin{equation}\label{e3.10}
      ||u^a||_{C^{2,\alpha}(\Omega^a_{1/8})}\leq C_{1/8}, \ \ \forall x\in\Omega^a_{1/8}
    \end{equation}
 holds for some $\alpha\in(0,1)$, where $C_{1/8}$ is a positive constant depending on $n,k,\gamma$ but not on $u^a,\Omega^a$.
\end{prop}

Using the relation between reverse isoperimetric inequality and $\gamma-$ball condition in Theorem \ref{t5.1}, there exist a sequence of $t=t_j\to\infty$ and a sequence of $t^*_j\in(t_j/2,t_j/3)$, such that the domain $\Omega_{t^*_j}$ satisfies uniformly $\gamma^*-$ball condition for some constant $\gamma^*$ independent of $j$. Applying uniform $\gamma^*$-ball condition and \eqref{e3.10} to $u^a$ yields that
   $$
    [u]_{C^{2,\alpha}(\Omega_{\sqrt{t_j}/8})}=t_j^{-\alpha/2}[u^a]_{C^{2,\alpha}(\Omega^a_{1/8})}\leq Ct_j^{-\alpha/2}.
   $$
After sending $t_j$ to infinity, we thus conclude that $u$ must be a quadratic function. Conversely, if $u(x)=a^{ij}x_ix_j$ is a quadratic polynomial for positive definite matrix $A=[a^{ij}]$, we have
   \begin{eqnarray*}
     {\mathcal{I}}(u)&=&\liminf_{t\to\infty}\frac{\int_{\Omega_t}\sqrt{a^{ik}a^{jk}x_ix_j}}{\Big(\int_{\Omega_t}|t-a^{ij}x_ix_j|^{\frac{n}{n-1}}\Big)^{\frac{n-1}{n}}}\\
     &\leq&C\liminf_{t\to\infty}\frac{t^{\frac{n+1}{2}}}{\big(t^{\frac{n}{2}+\frac{n}{n-1}}\big)^{\frac{n-1}{n}}}<\infty.
   \end{eqnarray*}
The necessarity and sufficiency in Theorem \ref{t1.3} have been shown. $\Box$\\

\vspace{40pt}

\section{Reduction of Hessian quotient equation by Legendre transformation}

The main obstacle in proving of Theorem \ref{t1.5} is the lacking of $C^2$-estimation for Hessian quotient equation \eqref{e1.7}. Fortunately, one can utilize Legendre transformation to change the Hessian quotient equation into $k-$Hessian equation \eqref{e1.10}. Then, with the help of $C^2-$estimation of \eqref{e1.10}, one derives an analogue Bernstein property for \eqref{e1.7}.

 Suppose that $u$ is a strict convex solution to \eqref{e1.7} on $\gamma-$ball domain $\Omega$, satisfying
   \begin{equation}\label{e4.1}
     \begin{cases}
       S_{n,l}(D^2u)=1, & \forall x\in\Omega\\
       u(x)=1, & \forall x\in\partial\Omega
     \end{cases}
   \end{equation}
and
   \begin{equation}\label{e4.2}
     u(0)=0, \ \ Du(0)=0.
   \end{equation}
 Denoting
    $$
     \Omega_{1/2}\equiv\Big\{x\in\Omega\Big|\ u(x)<1/2\Big\}
    $$
 to be the sub-level set of $u$ and
    $$
     v(y)\equiv\sup_{x\in\Omega_{1/2}}(x\cdot y-u(x)), \ \ \forall y\in\Omega^*\equiv Du(\Omega_{1/2})
    $$
 to be the Legendre transformation of $u$, there hold
    \begin{equation}\label{e4.3}
      y=Du(x), \ \ x=Dv(y),\ \ D^2u(x)=\Big(D^2v(y)\Big)^{-1}
    \end{equation}
 and
    \begin{equation}\label{e4.4}
      v(0)=0, \ \ Dv(0)=0.
    \end{equation}
So, $v$ must be a solution to
   \begin{equation}\label{e4.5}
     S_{n-l}(D^2v)=1, \ \ \forall y\in\Omega^*.
   \end{equation}
The next lemma shows that $\Omega$ is not large and $\Omega^*$ is not small.

\begin{lemm}\label{l4.1}
   Under the assumption of Theorem \ref{t1.5} and let $u$ be a strict convex solution of \eqref{e4.1}-\eqref{e4.2} after Pogorelov's first normalization. There exist positive constants $C_1$ and $C_2$ depending only on $n,l$ and $\gamma$, such that $\Omega$ is not large in the sense of
    \begin{equation}\label{e4.6}
      \{0\}\in\Omega\subset C_1B_1(x_0)
    \end{equation}
and $\Omega^*$ is not small in the sense of
    \begin{equation}\label{e4.7}
      \Omega^*\supset C_2^{-1}B_1(0).
    \end{equation}
\end{lemm}

\noindent\textbf{Proof.} By Newton-Maclaurin's inequality, it is inferred from \eqref{e4.1} that
   $$
     \det(D^2u)=S_l(D^2u)\geq C_n^l{\det}^{\frac{l}{n}}(D^2u)\Rightarrow\det(D^2u)\geq\Big(C_n^l\Big)^{\frac{n}{n-l}}.
   $$
Therefore, one concludes from Lemma \ref{l2.2} that the normalized radius satisfies
   \begin{equation}\label{e4.8}
     R\leq C
   \end{equation}
for some positive constant $C$ depending only on $n,l$ and $\gamma$. So, \eqref{e4.6} follows. Next, we show that $dist(0,\partial\Omega^*)$ is uniformly bounded from below. Take a direction $\omega\in{\mathbb{S}}^{n-1}$ such that
    $$
     r_\omega=\inf_{\omega'\in{\mathbb{S}}^{n-1}}r_{\omega'},\ \  r_{\omega'}\equiv\sup\Big\{r>0|\ r\omega'\in\Omega^*\Big\}.
    $$
Setting $y_\omega\equiv r_\omega\omega\in\partial\Omega^*$ to be the boundary point of $\Omega^*$ in direction $\omega$, there must be a boundary point $x_\omega\in\partial\Omega_{1/2}$ such that
   $$
    |Du|(x_\omega)=r_\omega.
   $$
Drawing a straight segment
   $$
    \beta=\beta(t), \ \ \beta(0)=0, \ \ \beta(t_\omega)=x_\omega,\ \  t\in[0,t_\omega]
   $$
parameterized by arc-length parameter $t$, one has
   $$
    t_\omega\leq C_1
   $$
by \eqref{e4.6}. Since $Du(\beta(t))\beta'(t)$ is monotone increasing due to the convexity of $u$,
   $$
    \frac{1}{2}=u(x_\omega)=\int^{t_\omega}_0Du(\beta(t))\beta'(t)dt\leq  r_\omega\cdot t_\omega\leq C_1r_\omega.
   $$
Hence, \eqref{e4.7} follows by setting $C_2\equiv 2C_1$ and the proof of the lemma was done. $\Box$\\

Next lemma gives $C^2-$bound of $v$ upon the assumptions of Theorem \ref{t1.5}.

\begin{lemm}\label{l4.2}
  Under the assumption of Theorem \ref{t1.5} and letting $u$ be a strict convex solution of \eqref{e4.1}-\eqref{e4.2} after Pogorelov's first normalization, there exists a positive constant $C_4$ such that
     \begin{equation}\label{e4.9}
        |D^2v(y)|\leq C_4, \ \ \forall y\in (2C_2)^{-1}B_1(0)
     \end{equation}
  holds for $C_2$ coming from Lemma \ref{l4.1}.
\end{lemm}

\noindent\textbf{Proof.} Rescaling the solution $v$ of \eqref{e1.16} by $w(x)=(R/2)^2v(x/(R/2))$ for $R=C_2^{-1}$ in Lemma \ref{l4.1} and shifting the center of ball, one actually derives that for any solution $w$ of \eqref{e1.16} on $B_R$, there holds
    \begin{equation}\label{e4.10}
      \sup_{B_{R/2}}|D^2w|\leq\psi(n,l,\sup_{B_{R}}(|Dw|/R)).
    \end{equation}
This thus implies the uniform boundedness of $|D^2v|$ in $B_{1/(2C_2)}$ by Lemma \ref{l4.1}. $\Box$\\

Combining Lemma \ref{l4.2} with Krylov's interior estimation \cite{K}, we have achieved the following {\it a-priori} bound for higher derivatives.

\begin{lemm}\label{l4.3}
Under the assumption of Theorem \ref{t1.5} and letting $u$ be a strict convex solution of \eqref{e4.1}-\eqref{e4.2} after Pogorelov's first normalization, there exists a positive constant $C_5$ such that
  \begin{equation}\label{e4.11}
     |D^3v(y)|\leq C_5, \ \ \forall y\in (4C_2)^{-1}B_1(0).
  \end{equation}
\end{lemm}

\noindent\textbf{Completing the proof of Theorem \ref{t1.5}} By Lemma \ref{l4.3}, for each indices $k, p, q=1,\cdots,n$,
   $$
     \sup_{B_{1/(4C_2)}}|D_kv_{pq}|\leq C, \ \ \forall k,p,q=1, 2,\cdots,n.
   $$
Using the relation
   \begin{eqnarray*}
    D_kv_{pq}&\equiv&\frac{\partial}{\partial y_k} v_{pq}(y)=\frac{\partial}{\partial x_l}u^{pq}(x)\frac{\partial x_l}{\partial y_k}\\
    &=&D_lu^{pq}\Bigg[\frac{\partial y_k}{\partial x_l}\Bigg]^{-1}=u^{kl}D_lu^{pq}, \ \ [u^{pq}]=[u_{pq}]^{-1}
   \end{eqnarray*}
for $u=u^a$ after first normalization, one gets that
   \begin{equation}\label{e4.12}
     |(u^a)^{kl}D_l(u^a)^{pq}|(0)\leq C, \ \ \forall k,p,q=1,2,\cdots,n.
   \end{equation}
After scaling $u^a$ back to $u$, it yields from \eqref{e4.12} and Theorem \ref{t5.1} that
   $$
    |u^{kl}D_lu^{pq}|(0)\leq Ct_j^{-1/2}
   $$
holds for some positive constant $C$ independent of $j$. Sending $j$ to infinity yields that
   $$
    u^{kl}D_lu^{pq}(0)=0, \ \ \forall k,p,q=1,2,\cdots, n.
   $$
Noting that the reverse isoperimetric inequality is invariant under subtraction by a tangential linear function
   $$
     l_{x_0}(x)\equiv u(x_0)+Du(x_0)(x-x_0)
   $$
for each $x_0\in {\mathbb{R}}^n$ in Proposition \ref{p5.1}, there holds
   \begin{equation}\label{e4.13}
    u^{kl}D_lu^{pq}(x_0)=0, \ \ \forall k,p,q=1,2,\cdots, n
   \end{equation}
for each $x_0\in {\mathbb{R}}^n$, which in turn implies that $u$ must be a quadratic polynomial. The proof of Theorem \ref{t1.5} was done. $\Box$\\

As shown above, a key ingredient in proof of Theorem \ref{t1.5} is an interior $C^2-$bound of the $(n-l)$-Hessian equation \eqref{e4.5}. Since the domain $\Omega^*$ may not contain any sub-level set of $v$, the theorem of Chou-Wang \cite{CW} under strict convexity can not be applied directly. Fortunately, with the help of Theorem \ref{t1.6}, one gets Theorem \ref{t1.7} as a corollary of Theorem \ref{t1.5}.

\vspace{40pt}

\section{Reverse isoperimetric inequality and $\gamma-$ball condition}

In this section, we will firstly prove the following relation between reverse isoperimetric inequality and $\gamma-$ball condition.

\begin{theo}\label{t5.1}
  Supposing that a locally strict convex function $u$ satisfies reverse isoperimetric inequality
    \begin{equation}\label{e5.1}
       \int_{\Omega_t}|D u|\leq \gamma\Bigg(\int_{\Omega_t}|t-u|^{\frac{n}{n-1}}\Bigg)^{\frac{n-1}{n}}, \ \ \Omega_t\equiv\Big\{x\in{\mathbb{R}}^n\Big|\ u(x)<t\Big\}
    \end{equation}
  for some constants $\gamma>0$ and $t$, then there exists $t'\in(t/3,t/2)$ such that $\Omega_{t'}$ satisfies $\gamma'-$ball condition for some positive constant $\gamma'$ depending only on $n$ and $\gamma$.
\end{theo}

\noindent\textbf{Remark.} The reverse isoperimetric inequality
     \begin{equation}\label{e5.2}
      {\mathcal{I}}(u)\equiv\liminf_{t\to\infty}\frac{\int_{\Omega_t}|D u|}{\Big(\int_{\Omega_t}|t-u|^{\frac{n}{n-1}}\Big)^{\frac{n-1}{n}}}<\infty
    \end{equation}
 is clearly invariant under Pogorelov's first normalization, subtracting or multiplying by constants. Since the sub-level set may not be bounded for general locally convex functions, the reverse isoperimetric inequality is not invariant under subtracting by arbitrary linear functions. Henceforth, we only consider reverse isoperimetric inequality for locally convex functions satisfying the normalized condition \eqref{e1.13}. However, we still have the following invariant.

 \begin{prop}\label{p5.1}
   Supposing that $u$ is a locally strict convex function on ${\mathbb{R}}^n$ satisfying reverse isoperimetric inequality \eqref{e5.2}, then the function
      $$
       v(x)\equiv u(x)-u(x_0)-Du(x_0)(x-x_0),\ \ \forall x\in{\mathbb{R}}^n
      $$
   also satisfies the reverse isoperimetric inequality \eqref{e5.2} for each $x_0\in {\mathbb{R}}^n$.
 \end{prop}

 Let us prove first an elementary lemma.

 \begin{lemm}\label{l5.1}
   Letting $u$ be a locally strict convex function satisfying
      \begin{equation}\label{e5.3}
       u(0)=0, \ \ Du(0)=0
      \end{equation}
and setting
   $$
    r\equiv|x|\in{\mathbb{R}}^+, \ \ \omega\equiv\frac{x}{|x|}\in{\mathbb{S}}^{n-1},
   $$
we have
      \begin{equation}\label{e5.4}
        \omega\cdot Du(r,\omega)\geq \inf_{\omega'\in{\mathbb{S}}^{n-1}}\omega'\cdot Du(1,\omega')\equiv\delta_1>0
      \end{equation}
and
     \begin{equation}\label{e5.5}
       u(r,\omega)\geq\delta_1(r-1)+\delta_2,\ \ \delta_2\equiv\inf_{\omega'\in{\mathbb{S}}^{n-1}}u(1,\omega')
     \end{equation}
 hold for all $x=r\omega\in B_1^c(0)$.
 \end{lemm}

 \noindent\textbf{Proof.} In terms of polar coordinates $(r,\omega)$ with $r\in{\mathbb{R}}^+, \omega\in{\mathbb{S}}^{n-1}$, we have
   \begin{eqnarray*}
     u_r&=&\omega\cdot  Du(r,\omega)\\
     u_{rr}&=&\omega\cdot D^2u(r,\omega)\cdot\omega>0.
   \end{eqnarray*}
 Therefore, one gets that
    $$
     u_r(r,\omega)\geq\delta_1, \ \ u(r,\omega)\geq\delta_1(r-1)+\delta_2
    $$
for each $r>1$ by monotonicity. The proof was done. $\Box$\\

Now, let us complete the proof of Proposition \ref{p5.1}. Fixing $x_0\in{\mathbb{R}}^n$, as one chooses a large $R>1$, there exists a small constant $\delta>0$ such that
   \begin{eqnarray}\nonumber\label{e5.6}
     &&u(x)\geq \delta|x|, \ \ \forall x\in B_R^c\\
     &&v(x)\geq\delta|x-x_0|\geq\frac{\delta}{2}|x|, \ \ \forall x\in B_R^c
   \end{eqnarray}
by convexities of $u$ and $v$. Hence,
   \begin{eqnarray}\nonumber\label{e5.7}
     \int_{\Omega^v_t}|Dv|&=&\Bigg(\int_{\Omega^v_t\cap B_R}+\int_{\Omega^v_t\cap B_R^c}\Bigg)|Du(x)-Du(x_0)|\\
      &\leq&C_{R,|u|_{C^1(B_R)}}+\Bigg(\frac{|Du(x_0)|}{\delta_1}+1\Bigg)\int_{\Omega^u_{tC_1}}|Du|
   \end{eqnarray}
holds by Lemma \ref{l5.1}, where
    \begin{equation}\label{e5.8}
      u(x)\leq v(x)+u(x_0)+|Du|(x_0)(|x|+|x_0|)\leq C_1t, \ \ \forall x\in \Omega^v_t\cap B_R^c
    \end{equation}
has been used for some positive constant $C_1$ depending only on $\delta, R$. Similarly, we claim that
   \begin{eqnarray}\nonumber\label{e5.9}
     \int_{\Omega^v_t}|t-v|^{\frac{n}{n-1}}&=&\Bigg(\int_{\Omega^v_t\cap B_R}+\int_{\Omega^v_t\cap B_R^c}\Bigg)\big|t-v(x)\big|^{\frac{n}{n-1}}\\
      &\geq&2^{-1}C_2^{\frac{n}{n-1}}\int_{\Omega^u_{t/C_2}}|t/C_2-u(x)|^{\frac{n}{n-1}}-C_{R,||u||_{C^1(B_R)}}
   \end{eqnarray}
for some positive constant $C_2$ and large $t$. In fact, a first observation is that
   \begin{equation}\label{e5.10}
      v(x)\leq u(x)+|Du|(x_0)(|x|+|x_0|)\leq C_2u(x), \ \ \forall x\in B^c_R
   \end{equation}
holds true for large $R$ and some positive constant $C_2$ to be chosen later so large that
   $$
    C_2\geq\frac{|Du|(x_0)+1}{\delta_1}+1,
   $$
where Lemma \ref{l5.1} has been used. Thus, there holds
   \begin{eqnarray*}
     \int_{\Omega^v_t}|t-v|^{\frac{n}{n-1}}&=&\Bigg(\int_{\Omega^v_t\cap B_R}+\int_{\Omega^v_t\cap B_R^c}\Bigg)\big|t-v(x)\big|^{\frac{n}{n-1}}\\
      &\geq&(1-\varepsilon)\int_{\Omega^v_t\cap B_R}t^{\frac{n}{n-1}}-C_{R,||u||_{C^1(B_R)}}+C_2^{\frac{n}{n-1}}\int_{\Omega^v_{t}\cap B_R^c}|t/C_2-u(x)|^{\frac{n}{n-1}}
   \end{eqnarray*}
by choosing $t$ large with respect to fixed $R$, where $\varepsilon>0$ is a small constant. When $x\in B_R^c$ and $u(x)\leq C_2^{-1}t$, it follows from \eqref{e5.6} that
  \begin{eqnarray*}
    v(x)&=& u(x)-u(x_0)-Du(x_0)(x-x_0)\\
     &\leq& C_2^{-1}t+|Du|(x_0)(|x|+|x_0|)\leq C_2^{-1}t+\frac{2|Du|(x_0)t}{C_2\delta}<t
  \end{eqnarray*}
by selecting $C_2$ large. So, we conclude that $\Omega^u_{t/C_2}\cap B_R^c\subset \Omega^v_t\cap B_R^c$. Another hand, for $x\in B_R$, $u(x)\leq C_2^{-1}t$ and $t$ large, one also has
  $$
    v(x)\leq \sup_{B_R}u+|Du|(x_0)(R+|x_0|)<t
  $$
and
  $$
   2^{-1}C_2^{\frac{n}{n-1}}\int_{\Omega^u_{t/C_2}\cap B_R}|t/C_2-u(x)|^{\frac{n}{n-1}}\leq\Bigg(\frac{1}{2}+\varepsilon\Bigg)\int_{\Omega^v_t\cap B_R}t^{\frac{n}{n-1}}+C_{R,||u||_{C^1(B_R)}}.
  $$
Hence, the claim \eqref{e5.9} was drawn. To proceed further, we need the following lemma.

\begin{lemm}\label{l5.2}
 There exists a positive constant $C_3$ depending only on $C_1,C_2$ and $n$, such that
    \begin{equation}\label{e5.11}
      \int_{\Omega^u_{t/C_2}}|t/C_2-u(x)|^{\frac{n}{n-1}}\geq C_3^{-1}\int_{\Omega^u_{tC_1}}|tC_1-u(x)|^{\frac{n}{n-1}}
    \end{equation}
 holds for convex function $u$ satisfying the normalized condition
    $$
     u(0)=0, \ \ Du(0)=0.
    $$
\end{lemm}

\noindent\textbf{Proof.} Using $\Omega^u_{t/C_2}$ as the base and origin as the vertex, we draw a cone whose boundary is just the graph of a function $w_1$ and denote $w_2\equiv (t/C_2-w_1(x))^{\frac{n}{n-1}}$ for short. By Fubini's theorem,
  \begin{eqnarray}\nonumber\label{e5.12}
   \int_{\Omega^u_{t/C_2}}|t/C_2-w_1(x)|^{\frac{n}{n-1}}&=&\int_{\Omega^u_{t/C_2}}\int^{(t/C_2)^{\frac{n}{n-1}}}_0\chi_{\{w_2(x)>s\}}(x,s)dsdx\\
    &=&\int^{(t/C_2)^{\frac{n}{n-1}}}_0|\chi_{\{w_2(x)>s\}}(x,s)|ds\\ \nonumber
    &=&\int^{(t/C_2)^{\frac{n}{n-1}}}_0|\Omega^u_{t/C_2}|\frac{(t/C_2-s^{\frac{n-1}{n}})^n}{(t/C_2)^n}ds\\ \nonumber
    &\geq& C^{-1}(t/C_2)^{\frac{n}{n-1}}|\Omega^u_{t/C_2}|
  \end{eqnarray}
holds for some positive constant $C$ depending only on $C_2$ and $n$. Using the convexity of $u$, we have also
   \begin{equation}\label{e5.13}
     |\Omega^u_{t/C_2}|\geq (C_1C_2)^{-n}|\Omega^u_{tC_1}|.
   \end{equation}
As a result, it follows from \eqref{e5.12} and \eqref{e5.13} that
  \begin{eqnarray*}
    \int_{\Omega^u_{t/C_2}}|t/C_2-u(x)|^{\frac{n}{n-1}}&\geq& C^{-1}(tC_1)^{\frac{n}{n-1}}|\Omega^u_{tC_1}|\\
     &\geq& C^{-1}\int_{\Omega^u_{tC_1}}|tC_1-u(x)|^{\frac{n}{n-1}}.
  \end{eqnarray*}
The proof of the lemma was done. $\Box$\\

Taking a sequence $t_j\to\infty$ such that
  $$
   \frac{\int_{\Omega_{t_j}}|Du|}{\Big(\int_{\Omega_{t_j}}|t-u|^{\frac{n}{n-1}}\Big)^{\frac{n-1}{n}}}\to {\mathcal{I}}(u)<\infty,
  $$
a combination of \eqref{e5.7}, \eqref{e5.9} with Lemma \ref{l5.2} yields that
    \begin{eqnarray*}
     \int_{\Omega^v_{t_j}}|Dv|\leq C\Bigg(\int_{\Omega^v_{t_j}}|{t_j}-v|^{\frac{n}{n-1}}+C_{R,||u||_{C^1(B_R)}}\Bigg)^{\frac{n-1}{n}}
    \end{eqnarray*}
for each $x_0\in {\mathbb{R}}^n$. Using the fact that
    $$
     \int_{\Omega^v_{t_j}}|{t_j}-v|^{\frac{n}{n-1}}\geq C^{-1}{t_j}^{\frac{n}{n-1}}|\Omega^v_{t_j}|
    $$
for large $j$, one obtains that ${\mathcal{I}}(v)<\infty$. So,
the proof of Proposition \ref{p5.1} was completed. $\Box$\\

We are now in a position to finish the proof of Theorem \ref{t5.1}. By performing the Pogorelov's first normalization, one needs only to consider the case $t=1$. For each $t\in(0,1)$, we set
    $$
     \mu(t)\equiv|\Omega_t|, \ \ \nu(t)\equiv|\partial\Omega_t|
    $$
for simplicity. By Fubini's theorem and Co-area formula, one has
   \begin{eqnarray}\nonumber\label{e5.14}
     \int_{\Omega_1}|1-u|^{\frac{n}{n-1}}&=&\frac{n}{n-1}\int_{\Omega_1}\int^1_{u(x)}(1-s)^{\frac{1}{n-1}}dsdx\\
     &=&\frac{n}{n-1}\int_{\Omega_1}\int^1_0(1-s)^{\frac{1}{n-1}}\chi_{\{u(x)<s\}}dsdx\\ \nonumber
     &=&\frac{n}{n-1}\int^1_0(1-s)^{\frac{1}{n-1}}\mu(s)ds
   \end{eqnarray}
and
   \begin{equation}\label{e5.15}
     \int_{\Omega_1}|D u|=\int^1_0\nu(s)ds.
   \end{equation}
Therefore, it follows from \eqref{e5.1} and \eqref{e5.14}-\eqref{e5.15} that
   \begin{equation}\label{e5.16}
     \int^1_0\nu(s)ds\leq C\gamma\Bigg(\int^1_0(1-s)^{\frac{1}{n-1}}\mu(s)ds\Bigg)^{\frac{n-1}{n}}.
   \end{equation}
By mean value theorem, there exists $s_*\in(1/3,1/2)$ such that
    \begin{equation}\label{e5.17}
       \frac{1}{6}\nu(s_*)=\int^{1/2}_{1/3}\nu(s)ds.
    \end{equation}
Since $u$ is a convex function, one has also
   \begin{equation}\label{e5.18}
     \mu(1)\geq\mu(s)\geq s^n\mu(1), \ \ \forall s\in(0,1).
   \end{equation}
Actually, to show that $\mu(s)\geq s^n\mu(1)$, we need only drawing a cone $V$ with vertex at $(0,0)\in{\mathbb{R}}^{n+1}$ and base $\Omega_1\times\{1\}$. By convexity of $u$, the section
  $$
   {\mathcal{S}}_s\equiv\Big\{x\in{\mathbb{R}}^n\Big|\ (x,s)\in V\Big\}
  $$
is contained inside $\Omega_s$. Moreover, its area is exactly given by $s^n\mu(1)$ by similarity. We finally arrive at the inequality $\mu(s)\geq s^n\mu(1)$. Summing \eqref{e5.16}-\eqref{e5.18} yields that
    \begin{equation}\label{e5.19}
       \nu(s_*)\leq C\gamma\mu^{\frac{n-1}{n}}(s_*).
    \end{equation}
So, the conclusion of Theorem \ref{t5.1} follows from the following lemma.

\begin{lemm}\label{l5.3}
  Under the assumptions of Theorem \ref{t5.1}, if \eqref{e5.19} holds for some positive constant $C$, there exists another positive constant $C'$ depending only on $n,k$ and $C$ such that $\Omega_{s_*}$ is a $\gamma'-$ball domain for $\gamma'=C'\gamma^{n/2}$.
\end{lemm}

\noindent\textbf{Proof.} By John's lemma, we may assume that $\Omega_{s_*}$ satisfies the John's ball condition \eqref{e2.2} after some affine transformation $A, \det(A)=1$ with eigenvalues
   $$
     \mu(A)=(\mu_1,\mu_2,\cdots,\mu_n), \ \ 0<\mu_1\leq\mu_2\leq\cdots\leq\mu_n.
   $$
To show the conclusion of Lemma \ref{l5.3}, one needs only to estimate the quotient of $\mu_n/\mu_1$ from above. Noting that for some universal constants $C_{n,1}$ and $C_{n,2}$, the volume of ellipsoid $A^{-1}(B_R)$ is given by
    \begin{equation}\label{e5.20}
       V(A^{-1}(B_R))=C_{n,1}S_n(\mu(A^{-1}))R^n
    \end{equation}
and the surface area of ellipsoid $A^{-1}(B_R)$ satisfies that
    \begin{equation}\label{e5.21}
      C_{n,2}^{-1}S_{n-1}(\mu(A^{-1}))R^{n-1}\leq S(A^{-1}(B_R))\leq C_{n,2}S_{n-1}(\mu(A^{-1}))R^{n-1}.
    \end{equation}
By volume and area comparisons of convex bodies, it follows from \eqref{e5.19}-\eqref{e5.21} that
    \begin{eqnarray}\nonumber\label{e5.22}
        &&S_{n-1}(\mu(A^{-1}))\leq C\gamma[S_n(\mu(A^{-1}))]^{\frac{n-1}{n}}\\
        &\Leftrightarrow& [S_1(\mu(A))]^n\leq C\gamma^nS_n(\mu(A)).
    \end{eqnarray}
Hence, one obtains
     \begin{equation}\label{e5.23}
       \mu_n\leq C\gamma^n\mu_1
     \end{equation}
for some positive constant $C$. The proof was done. $\Box$\\

\vspace{40pt}

\section{Volume growth and uniformly $\gamma-$ball condition}

Key ingredient in proving of Theorem \ref{t1.3}, \ref{t1.5} and \ref{t1.7} is the verification of uniform $\gamma-$ball condition using the reverse isoperimetric inequality. In this section, we will prove the $\gamma-$ball condition by an {\it a-priori} bound on the volume of the sublevel set.

\begin{theo}\label{t6.1}
 Suppose that $\Omega\subset{\mathbb{R}}^n$ is a bounded convex domain with volume $|\Omega|$, and $u\in C^2(\Omega)$ is a convex solution to $k-$Hessian equation
    \begin{equation}\label{e6.1}
     \begin{cases}
      S_k(D^2u)=1, & \forall x\in\Omega\\
        u(0)=0, Du(0)=0, u(x)=1, & \forall x\in\partial\Omega
     \end{cases}
    \end{equation}
 for some $1\leq k\leq n-1$. Then $\Omega$ satisfies the $\gamma-$ball condition for some positive constant $\gamma=\gamma(n,k,|\Omega|)$ depending only on $n,k$ and $|\Omega|$.
\end{theo}

\noindent\textbf{Proof.} By John's lemma, there exists a matrix $A\in SL(n)$ with eigenvalues
  $$
    \mu(A)=(\mu_1,\cdots,\mu_n),\ \ 0<\mu_1\leq\mu_2\leq\cdots\leq\mu_n, \ \ \Pi_{i=1}^n\mu_i=1
  $$
such that
    \begin{equation}\label{e6.2}
       C_n^{-1}E_{A,R}(x_0)\subset \Omega\subset C_nE_{A,R}(x_0)
    \end{equation}
holds for ellipsoid
    $$
     E_{A,R}(x_0)\equiv\Bigg\{x\in{\mathbb{R}}^n\Bigg|\ \frac{|x_1-x_{0,1}|^2}{\mu_1^2}+\cdots+\frac{|x_n-x_{0,n}|^2}{\mu_n^2}=R^2\Bigg\},
    $$
where $x_0=(x_{0,1},\cdots,x_{0,n})$ and
 $$
     C_n^{-1}|\Omega|\leq R^n\leq C_n|\Omega|.
 $$
Direct computation shows that the upper barrier function
   $$
    v(x)=\frac{1}{2S_k^{1/k}(\mu((A^{-1})^2))}\Bigg[\Bigg(\frac{|x_1-x_{0,1}|^2}{\mu_1^2}+\cdots+\frac{|x_n-x_{0,n}|^2}{\mu_n^2}\Bigg)-C_n^{-2}R^2\Bigg]+1,
   $$
satisfies that
   \begin{equation}\label{e6.3}
     \begin{cases}
       S_k(D^2v)=1, & \forall x\in C_n^{-1}E_{A,R}(x_0)\\
       v(x)=1, & \forall x\in\partial(C_n^{-1}E_{A,R}(x_0)).
     \end{cases}
   \end{equation}
Comparison of $u$ with $v$ yields that
   $$
    0\leq u(x)\leq v(x)\Rightarrow 0\leq v(x_0)=-\frac{C_n^{-2}R^2}{2S_k^{1/k}(\mu((A^{-1})^2))}+1.
   $$
One gets thus that
   \begin{equation}\label{e6.4}
      R^2\leq 2C_n^2S_k^{1/k}(\mu((A^{-1})^2)).
   \end{equation}
On another hand, after constructing a lower barrier function
   $$
    w(x)=\frac{1}{2S^{1/k}_k(\mu((A^{-1})^2))}\Bigg[\Bigg(\frac{|x_1-x_{0,1}|^2}{\mu_1^2}+\cdots+\frac{|x_n-x_{0,n}|^2}{\mu_n^2}\Bigg)-C_n^{2}R^2\Bigg]+1,
   $$
one has
   \begin{equation}\label{e6.5}
     \begin{cases}
       S_k(D^2w)=1, & \forall x\in \Omega\\
       w(x)\leq1, & \forall x\in\partial\Omega.
     \end{cases}
   \end{equation}
A similar comparison of $u$ with $w$ yields that
   $$
    0=u(0)\geq w(0)=-\frac{C_n^2R^2}{2S^{1/k}_k(\mu((A^{-1})^2))}+1
   $$
and so
   \begin{equation}\label{e6.6}
      R^2\geq 2C_n^{-2}S^{1/k}_k(\mu((A^{-1})^2)).
   \end{equation}
Noting that by \eqref{e6.6} and $k\leq n-1$,
   $$
    S_{n-1}^{\frac{1}{n-1}}(\mu((A^{-1})^2))\leq CS_k^{\frac{1}{k}}(\mu((A^{-1})^2))\leq CR^2.
   $$
Together with $S_n(\mu((A^{-1})^2))=1$, one concludes that
   $$
    \mu_n^{-2}\geq C^{-1}R^{2-2n}
   $$
and thus
    \begin{equation}\label{e6.7}
      \mu_n\leq CR^{n-1}, \ \ \mu_1\geq C^{-1}R^{-(n-1)^2}.
    \end{equation}
As mentioned in Section 2, $\Omega$ satisfies the $\gamma-$ball condition for
   $$
    \gamma\leq CR^{n(n-1)/2}\leq C|\Omega|^{(n-1)/2}.
   $$
This is exactly the desired conclusion of Theorem \ref{t6.1}. $\Box$\\

In order completing the proof of Theorem \ref{t1.8}, one needs also another version of Theorem \ref{t6.1} for Hessian quotient equation as following.

\begin{theo}\label{t6.2}
 Suppose that $\Omega\subset{\mathbb{R}}^n$ is a bounded convex domain with volume $|\Omega|$, and $u\in C^2(\Omega)$ is a convex solution to Hessian quotient equation
    \begin{equation}\label{e6.8}
     \begin{cases}
      S_{n,l}(D^2u)=1, & \forall x\in\Omega\\
        u(0)=0, Du(0)=0, u(x)=1, & \forall x\in\partial\Omega
     \end{cases}
    \end{equation}
 for some $1\leq l\leq n-1$. Then $\Omega$ satisfies the $\gamma-$ball condition for some positive constant $\gamma=\gamma(n,l,|\Omega|)$ depending only on $n,l$ and $|\Omega|$.
\end{theo}

Theorem \ref{t6.2} can be proven as that in Theorem \ref{t6.1} without difficulty by constructing upper and lower barrier functions. Now, let us complete the proof of Theorem \ref{t1.8} as follows. Upon the assumption \eqref{e1.20}, there exists a sequence of $t_j\to\infty$, such that
    $$
     t_j^{-n/2}|\Omega_{t_j}|\leq C<\infty, \ \ \forall j.
    $$
Performing Pogorelov's first normalization $u_j^a(x)=t_j^{-1}u(\sqrt{t_j}x)$, one gets a sequence of sub-level sets $\Omega^a_{j,1}$ of $u_j^a$, which have uniformly bounded volumes. Utilizing Theorem \ref{t6.1} and \ref{t6.2}, we conclude a uniformly $\gamma-$ball condition for these domains. Repeating the arguments in Theorem \ref{t1.3}, \ref{t1.5} and \ref{t1.7} gives the Bernstein property of Hessian equations \eqref{e1.7} or \eqref{e1.10}, which thus implies the reverse isoperimetric inequality \eqref{e1.14}. Vice versa, if the reverse isoperimetric inequality \eqref{e1.14} holds, then the solution must be a quadratic function. So \eqref{e1.20} is clearly true. We have now completed the proof of Theorem \ref{t1.8} upon the condition \eqref{e1.20}. The validity of Theorem \ref{t1.8} upon the condition \eqref{e1.21} is also not hard to be verified using the Fubini's theorem
   \begin{equation}\label{e6.9}
     t^{-p-n/2}\int_{\Omega_{t-1}}|t-u|^p=pt^{-p-n/2}\int^{t}_0(t-s)^{p-1}|\Omega_s\cap\Omega_{t-1}|ds
   \end{equation}
and mean value theorem for integral. The proof of Theorem \ref{t1.8} was done. $\Box$\\

\vspace{40pt}

\section*{Acknowledgments}

The author would like to express his deepest gratitude to Professors Xi-Ping Zhu, Kai-Seng Chou, Xu-Jia Wang and Neil Trudinger for their constant encouragements and warm-hearted helps. This paper was also dedicated to the memory of Professor Dong-Gao Deng.

\end{document}